\documentclass[11pt,twoside,reqno]{amsart}

\usepackage{amssymb}

 \newtheorem{thm}{Theorem}[section]

 \theoremstyle{definition}
 \newtheorem{defn}[thm]{Definition}
 \theoremstyle{remark}
 \newtheorem{rem}[thm]{Remark}
 
 \numberwithin{equation}{section}

\def\R{{\mathbb R}}

\def\Om{\Omega}

\def\al{\alpha}
\def\be{\beta}
\def\sbs{\subset}
\def\vep{\varepsilon}

\begin{document}
\setcounter{page}{1}

%Topmatter

\title[]{A note on $H$-convergence}

\author[]{Bj\"{o}rn Gustafsson and Jacqueline Mossino}

\address{%%
Mathematics Department,
Royal Institute of Technology,
S-10044 Stockholm,
Sweden.}

\email{gbjorn@kth.se}

\address{C.M.L.A., Ecole Normale Sup\'erieure de Cachan,  
61, Avenue du Pr\'esident Wilson, 
94235 Cachan cedex, France. }

\email{Jacqueline.Mossino@cmla.ens-cachan.fr}

\thanks{The first author is supported by the Swedish Research Council}



\date{August 10, 2006}

\begin{abstract}
We give a criterion for $H$-convergence of conductivity
matrices in terms of ordinary weak convergence of the factors
in certain quotient representations of the matrices.
\end{abstract}

\maketitle

\section{Introduction}

Questions of homogenization of rapidly varying coefficients
in elliptic partial differential equations have been considered
by mathematicians at least since the 1970s, and by physicists
and engineers much longer (cf. the references in
\cite{Spagnolo 68}, \cite{Bensoussans-Lions-Papanicolaou 78},
\cite{Sanchez-Palenzia 80}). A typical example is when the
conductivity matrix $A^\vep=A^\vep (x)$ in an equation
\begin{equation}\label{de}
-{\rm div\,} (A^\vep \nabla u^\vep) = f
\end{equation}
in some (bounded) domain $\Omega\subset \R^n$ 
oscillates rapidly at a length scale $\vep >0$ and one wants to identify
a limiting matrix $A=A(x)$ (presumably less oscillating than the
$A^\vep$) such that, as $\vep\to 0$,
the solutions $u^\vep$ converge in some weak sense
to the solution $u$ of the corresponding homogenized
equation:
$$
-{\rm div\,} (A \nabla u) = f.
$$
 
In the 1970s and 80s, F.~Murat, L.~Tartar identified the 
appropriate type of convergence, $H$-convergence, for the above
type of problems
and started developing general theories for it. 
Earlier work was much concerned with
special cases, like strictly periodic structures. 
In this little note we give an equivalent condition for general
$H$-convergence in terms of ordinary weak convergence
for the factors in certain quotient
representations of the conductivity matrices, namely for matrices
$M^\vep$ and $P^\vep$ appearing when writing $M^\vep A^\vep =P^\vep$.
In the special case of stratified media ($A^\vep$ depending on only one
of the coordinates) and certain generalizations thereof, 
explicit decompositions
of this type have been constructed and used for proving
$H$-convergence in a series of papers \cite{Dufour-Fabre-Mossino 96},
\cite{Fabre-Mossino 98}, \cite{Courilleau 01}, 
\cite{Gustafsson-Mossino 03}. The purpose of this note is
to point out that the existence of such quotient
representations is a completely general fact in connection
with $H$-convergence $A^\vep\to A$. 

The proof consists of an adaptation of methods developed by
F.~Murat and L.~Tartar, e.g., in \cite{Murat 77}, \cite{Murat 78},
\cite{Tartar 79}, \cite{Tartar 85}, \cite{Tartar 86}, 
\cite{Murat-Tartar 97}. In fact, even the result can be said to be
implicit in their work, but perhaps not explicit.

For an extension of the results in this note to $H$-convergence in linear elasticity, 
see \cite{Gustafsson-Mossino 06}.

%%%%%%%%%%%%%%%%%%%%%%%%%%%%%%%%%%%%%%%%%%%%%%%%%%%%%%%

\section{$H$-convergence}

\begin{defn}
Let $0<\al \leq\be <\infty$, $\Om\sbs \R^n$ (a bounded domain). 
Then $M(\al,\be;\Om)$
denotes the set of invertible real-valued $n\times n$ matrices $A=A(x)$
with entries in $L^\infty (\Om)$ and satisfying almost 
everywhere in $\Om$ the estimates
$$
(A\xi,\xi)\geq \al |\xi|^2,          
$$
$$
(A^{-1}\xi,\xi)\geq \be^{-1} |\xi|^2       
$$
for $\xi \in \R^n$. 
\end{defn}

Above the bracket $(\cdot ,\cdot )$ denotes the scalar product
in $\R^n$: $(\xi, \eta) =\sum_j \xi_j \eta_j$.
Divergence and curl of vector fields are defined as usual:
${\rm div\,}D$ is the scalar
$$
{\rm div\,}D = \sum_j \frac{\partial D_j}{\partial x_j},
$$
and ${\rm curl\,} E$ is the antisymmetric tensor with components
$$
({\rm curl\,} E )_{ij} = \frac{\partial E_i}{\partial x_j}
-\frac{\partial E_j}{\partial x_i}.
$$
We shall also need to take ${\rm div\,}$ and ${\rm curl\,}$
of matrices, and then the above definitions apply
to the row vectors, i.e., to  the
last index. Thus, with $M= (M_{ij})$, $P=(P_{ij})$,
$$
({\rm div\,}P)_i = \sum_j \frac{\partial P_{ij}}{\partial x_j}
$$
$$
({\rm curl\,} M )_{ijk} = \frac{\partial M_{ij}}{\partial x_k}
-\frac{\partial M_{ik}}{\partial x_j}.
$$

The parameter $\vep>0$ to be used from now on
is by convention restricted to take values 
only in a sequence tending to zero
(e.g., $\vep\in \{ 1,\frac{1}{2}, \frac{1}{3},\dots \}$).

\begin{defn}
(Tartar \cite{Tartar 85})
Let $A, A^\vep \in M(\al,\be;\Om)$ for some $0<\al \leq\be <\infty$
and all $\vep >0$. 
Then $A^\vep$ is said to  $H$-converge to $A$,
$$
A^\vep \stackrel{H}{\longrightarrow} A,
$$
as $\vep\to 0$ if the following holds.
Whenever vector fields $D^\vep$, $D$, $E^\vep$, $E \in L^2(\Om)^n$
satisfy
\begin{equation}\label{DAEvep}
D^\vep  =A^\vep E^\vep,
\end{equation}
\begin{equation}\label{DvepD}
D^\vep \rightharpoonup D \text{ weakly in }L^2 (\Om)^n,
\end{equation}
\begin{equation}\label{EvepE}
E^\vep \rightharpoonup E \text{ weakly in }L^2 (\Om)^n
\end{equation}
with
\begin{equation}\label{divDcomp}
\{{\rm div \,} D^\vep\}_{\vep >0}
\text{  relatively compact in } H^{-1}(\Om),
\end{equation}
\begin{equation}\label{curlEcomp}
\{{\rm curl \,} E^\vep\}_{\vep >0}
\text{  relatively compact in } H^{-1}(\Om)^{n\times n}
\end{equation} 
then
\begin{equation}\label{DAE}
D=AE.
\end{equation}
\end{defn}

It is well-known \cite{Murat 77} that the $H$-limit is unique,
that the set $M(\al,\be;\Om)$ is sequentially
compact for $H$-convergence, and also that 
$H$-convergence is stable under transposition of the matrices: 
if $A^\vep \stackrel{H}{\longrightarrow} A$ then
$^t\! A^\vep \stackrel{H}{\longrightarrow} \,^t\!A$.

%%%%%%%%%%%%%%%%%%%%%%%%%%%%%%%%%%%%%%%%%%%%%%%%%%%%%%%%%%%%%%%%%%%

\section{The result} 

\begin{thm}
Let $A, A^\vep \in M(\al, \be; \Om)$ for some $0<\al\leq\be <\infty$.
Then $A^\vep \stackrel{H}{\longrightarrow} A$
as $\vep \to 0$ if and only if there exist
$n\times n$ matrices $M^\vep$, $M$, $P^\vep$, $P$ with entries in
$L^2 (\Om)$ and with $M$ (and hence $P$) invertible, such that 
\begin{equation}\label{MAPvep}
M^\vep A^\vep= P^\vep, 
\end{equation}
\begin{equation}\label{MAP}
MA=P,
\end{equation}
\begin{equation}\label{Mconv}
M^\vep \rightharpoonup M \text{ weakly in } L^2 (\Om)^{n\times n},
\end{equation}
\begin{equation}\label{Pconv}
P^\vep \rightharpoonup P \text{ weakly in } L^2 (\Om)^{n\times n}.
\end{equation}
with
\begin{equation}\label{Mcomp}
\{{\rm curl \,} M^\vep\}_{\vep >0}
\text{ relatively compact in } H^{-1}(\Om)^{n\times n\times n},
\end{equation} 
\begin{equation}\label{Pcomp}
\{{\rm div \,} P^\vep\}_{\vep >0}
\text{ relatively compact in } H^{-1}(\Om)^n.
\end{equation} 
When this is the case $M$ can be chosen 
to be the identity matrix I
and $M^\vep$ so that ${\rm curl \,} M^\vep =0$. 
%(Alternatively, $P$ can chosen arbitrarily and
%$P^\vep$ so that ${\rm div \,} P^\vep =0$.)
\end{thm}

\begin{proof}
The proof is based on the ``div-curl lemma'' of compensated
compactness \cite{Murat 77}, \cite{Tartar 85},
\cite{Evans 90}. We recall that this
lemma in general says that if
$f^\varepsilon, g^\varepsilon, f,g \in L^2(\Omega)^n$
are vector fields such that $f^\varepsilon \rightharpoonup f$,
$g^\varepsilon \rightharpoonup g$ weakly in $L^2(\Omega)^n$
and such that ${\rm div\,} f^\varepsilon$ and the components of  
${\rm curl\,} g^\varepsilon$
are all contained in a compact subset
of $H^{-1} (\Omega)$, then 
$(f^\varepsilon, g^\varepsilon) \rightharpoonup (f, g)$
weakly as distributions.

First we prove the ``if''-part of the theorem, which is very easy. 
So assume we have the decompositions
(\ref{MAPvep}) and (\ref{MAP}) with weak convergences 
$M^\vep\rightharpoonup M$ and
$P^\vep\rightharpoonup P$ as in the statement. We consider
$E^\vep$ and $D^\vep = A^\vep E^\vep$ satisfying 
(\ref{DvepD}), (\ref{EvepE}), (\ref{divDcomp}), (\ref{curlEcomp}).
Then (\ref{MAPvep}) acting on $E^\vep$ gives
$$
M^\vep D^\vep = P^\vep E^\vep,
$$
or, in components,
$$
\sum_j M^\vep_{ij} D^\vep_j = \sum_j P^\vep_{ij} E^\vep_j.
$$
Here the div-curl lemma applies for each $i$ and it follows 
that each of the members converge in the sense of distribution, to
$MD$ and $PE$ respectively. Thus we get
\begin{equation}\label{MDPE}
MD =PE,
\end{equation}
which, since $M$ is invertible, is the same as (\ref{DAE}).

Now we prove the ``only if'' part. First we have to construct 
the matrices $M^\vep, M, P^\vep, P$. We may take
$M=I$, $P=A$. Let $^t\!A$ denote the transpose of 
$A$ and let $e_i$ be the $i$:th unit column vector.  Thus
$e_i =  \nabla u_i$, where $u_i$ is the $i$:th coordinate function:
$$
u_i (x)= x_i.
$$
Setting also
$
f_i= {\rm div\,} (^t\! A e_i )
$
the equation
$$
{\rm div\,} (^t\! A \nabla u) = f_i
$$
is trivially solved by $u=u_i$.

Now, with $f_i$, $u_i$ as above ($1\leq i\leq n$),
there is for each $\vep >0$ a unique solution $u_i^\vep$ 
of the elliptic boundary value problem
\begin{equation}\label{BVP}
\begin{cases}
{\rm div\,} (^t\! A^\vep \nabla u_i^\vep) = f_i,\\
u_i^\vep - u_i \in H^1_0 (\Om).
\end{cases}
\end{equation}
Using it we define $M^\vep $ to be the matrix whose $i$:th
row is $^t\! \nabla u^\vep_i$. In other words, $M^\vep$ is the 
matrix with entries
$$
M^\vep_{ij} = \frac{\partial u^\vep_i}{\partial x_j}.
$$
Then we take $P^\vep$ to be
$$
P^\vep = M^\vep A^\vep,
$$
so that 
$^t\!P^\vep = ^t\!A^\vep (\nabla u_1^\vep,\dots, \nabla u_n^\vep)$.
The so defined matrices $M^\vep$, $P^\vep$ satisfy
$$
\begin{cases}
{\rm curl\,} M^\vep =0, \\
{\rm div\,} P^\vep = f,
\end{cases}
$$
where $f$ is the vector with components
$f_i\in H^{-1} (\Om)$.
In particular, the components of
${\rm curl\,} M^\vep$ and ${\rm div\,} P^\vep$
stay within a compact subset of $H^{-1}(\Om)$.

From (\ref{BVP}) we get for each $1\leq i\leq n$
the elliptic estimates
$$
\| u_i^\vep \|_{H^1(\Om)} \leq C <\infty,
$$
$$
\| ^t\!A^\vep\nabla u_i^\vep \|_{L^2(\Om)^n} \leq C <\infty.
$$
Thus for some subsequence of $\{\vep\}$ and some limit fields
$v_i$ and $\sigma_i$ we have convergences
\begin{equation}\label{convu}
u^\vep_i \rightharpoonup v_i, \text{ weakly in } H^1(\Om),
\end{equation}
\begin{equation}\label{convnablau}
\nabla u^\vep_i \rightharpoonup \nabla v_i, \text{ weakly in } L^2(\Om)^n,
\end{equation}
\begin{equation}\label{convsigma}
^t\!A^\vep \nabla u^\vep_i \rightharpoonup \sigma_i, 
\text{ weakly in } L^2(\Om)^n.
\end{equation}
The latter convergence together with (\ref{BVP}) shows that 
$$
{\rm div\,} \sigma_i = f_i.
$$

At this point we use the mentioned fact  
that $H$-convergence carries over to the transposed matricies.
Thus $^t\! A^\vep \stackrel{H}{\longrightarrow}\, ^t\!A$, and  
since ${\rm curl\,} \nabla u^\vep_i =0$
and ${\rm div\,} (^t\!A^\vep) = f_i$ are compact in $H^{-1} (\Om)$
it follows from the definition of this $H$-convergence that
$$
\sigma_i = ^t\!A \nabla v_i.
$$
Therefore $v_i$ solves the boundary value problem
\begin{equation}\label{BVPlimit}
\begin{cases}
{\rm div\,} (^t\! A \nabla v_i) = f_i,\\
v_i - u_i \in H^1_0 (\Om).
\end{cases}
\end{equation}

But this problem has the unique solution $u_i$. 
Thus we conclude that $v_i =u_i$ and that $\sigma_i =^t\!A\nabla u_i$.
It also follows that in (\ref{convu})--(\ref{convsigma}) 
we have convergence for the full sequence $\vep$
(because otherwise one
could extract a subsequence producing a different
solution of (\ref{BVPlimit})).
With this in mind, the convergences (\ref{convnablau}),
(\ref{convsigma}) state exactly that
$$
M^\vep \rightharpoonup M \text{ weakly in } L^2 (\Om)^{n\times n},
$$
$$
P^\vep \rightharpoonup P \text{ weakly in } L^2 (\Om)^{n\times n}.
$$
This proves the theorem.

\end{proof}

\begin{rem} 
It is clear that the matrices $M^\vep$, $M$, $P^\vep$, $P$
appearing in the decompositions (\ref{MAPvep}), (\ref{MAP}) are far
from being uniquely determined by $A^\vep$, $A$, even when
all the conditions (\ref{Mconv})--(\ref{Pcomp}) are satisfied.
For example, none of the conditions
(\ref{MAPvep})--(\ref{Pcomp}) are affected if
$M^\vep$, $M$, $P^\vep$, $P$
are multiplied from the left by one and the same invertible 
matrix $R=R(x)$ with  bounded Lipschitz coefficients, 

For the above reasons one cannot formulate the theorem as saying that
if (\ref{MAPvep}), (\ref{MAP}), (\ref{Mcomp}), (\ref{Pcomp}) hold
with $M$ invertible, then $A^\vep \stackrel{H}{\longrightarrow} A$
if and only if (\ref{Mconv}), (\ref{Pconv}) hold. 
However, such a statement is true if an appropriate normalization
is imposed.
Examples of such normalizations are that $M=I$ or that $P=I$.
One could also move one of the conclusions to be an assumption instead.
For example, the following statement is correct (and easy to deduce from
the theorem): if
(\ref{MAPvep}), (\ref{MAP}), (\ref{Mcomp}), (\ref{Pcomp}), 
(\ref{Mconv}) hold, then
$A^\vep \stackrel{H}{\longrightarrow} A$ if and only if
(\ref{Pconv}) holds. 

\end{rem}

\begin{rem}
From the formulation of the theorem one can easily pass to construction
of ``correctors'' (cf. \cite{Tartar 85}, \cite{Tartar 86}). 
Indeed, in order to construct correctors for $A^\vep$ one applies the
theorem to the transposed matrices: if 
$A^\vep \stackrel{H}{\longrightarrow} A$
then
$^t\!A^\vep \stackrel{H}{\longrightarrow}\, ^t\!A$.
Thus there are matrices $^t\!N^\vep$,
$^t\!N$, $^t\! Q^\vep$,  $^t\!Q$ such that 
$^t\!N^\vep \, ^t\!A^\vep = ^t\!Q^\vep$, 
$^t\!N \, ^t\!A = ^t\!Q$, 
$^t\!N^\vep \rightharpoonup ^t\! N$,    
$^t\!Q^\vep \rightharpoonup ^t\!Q$   with     
${\rm curl\,}^t\!N^\vep$ and ${\rm div\,}^t\!Q^\vep$
relatively compact in $H^{-1} (\Omega)$.
Here we choose the normalization $^t\! N=I$, $^t\! Q= A$.
If now we have vector fields $D^\vep$, $E^\vep$ as in the definition of 
the $H$-convergence (for $A^\vep$) and if an additional weak condition
is satisfied, e.g., that $N^\varepsilon$,
$Q^\varepsilon$ are bounded in $L^\infty(\Omega)^{n\times n}$, 
then the assertion is that
$$
E^\vep - N^\vep E \to 0 \quad \text{ strongly in } 
L_{\rm loc}^2 (\Omega)^n,
$$
$$
D^\vep - Q^\vep E \to 0 \quad \text{ strongly in } 
L_{\rm loc}^2 (\Omega)^n.
$$
This means that $N^\vep E$ and $Q^\vep E$ are good approximations
(correctors) of $E^\vep$ and $D^\vep$ respectively.

To prove the assertion, first notice that
$D^\vep - Q^\vep E = A^\vep (E^\vep - N^\vep E )$
and that therefore, for every $\omega \subset \Om$,
the $L^2 (\omega)^n$-norms of both $E^\vep - N^\vep E$
and $D^\vep - Q^\vep E$ can be estimated from above and below by
$\int_\omega (D^\vep - Q^\vep E , E^\vep - N^\vep E)\,dx$.
But the div-curl lemma gives that
$$
(D^\vep - Q^\vep E , E^\vep - N^\vep E)
=(D^\vep, E^\vep ) - (D^\vep,  N^\vep E)
-(Q^\vep E, E^\vep ) + (Q^\vep E,  N^\vep E)
$$
$$
\rightharpoonup (D, E) - (D, NE) -(QE, E) + (QE, NE) 
= (D-QE, E-NE) = 0
$$
in the sense of distributions (see \cite {Gustafsson-Mossino 06} for further
details).
From this the assertion follows.
\end{rem}

%%%%%%%%%%%%%%%%%%%%%%%%%%%%%%%%%%%%%%%%%%%%%%%%%%%

\section{Example}

The criterion of $H$-convergence in the theorem above is particularly 
useful in cases where it is possible to find the matrices
$M^\vep$ and $P^\vep$ {\it a priori} (without solving any Dirichlet
problem, e.g.). The main example for which this occurs is the 
case of stratified media, i.e., when $A^\vep$ depends on only one
of the coordinates, say $x_1$:
$$
A^\vep = A^\vep(x_1).
$$
Then the classical philosophy \cite{Murat 77}, \cite{Tartar 86}
is that  one should write the relation $D^\vep = A^\vep E^\vep$
in such a way that the ``bad'' components of $D^\vep$ and $E^\vep$
are expressed in terms of the ``good'' ones.
The good components are those for which one has control over the 
oscillations via the differential equation (\ref{de}) or via compactness
assumptions (\ref{divDcomp}), (\ref{curlEcomp}). 

In the stratified case, $D^\vep_1$ and $E^\vep_2$, \dots, $E^\vep_n$
are good and the rest are bad (see \cite{Tartar 86} for explanations),
and one thus writes (\ref{DAEvep}) as (suppressing $\vep$ for a moment)
$$
\begin{cases}
E_1 = \frac{1}{A_{11}} D_1 -\sum_{j\geq 2} \frac{A_{1j}}{A_{11}} E_j, \\
D_i =\frac{A_{i1}}{A_{11}} D_1 
+ \sum_{j\geq 2} (A_{ij} -\frac{A_{i1} A_{1j}}{A_{11}} )E_j \quad (i\geq 2).
\end{cases}
$$
In order to write this as $MD = PE$ one should take care to multiply
the bad quantities $E_1$ and $D_2$ \dots $D_n$ only by
good coefficients, for example constants. The simplest and
most natural choice of matrices $M$ and $P$ then is
$$
M = \left(\begin{array}{cc} \frac{1}{A_{11}} & 0 \\ \\ 
-\frac{A_{i1}}{A_{11}} & \delta_{ij} \end{array}\right), 
$$
$$
P =
\left(\begin{array}{cc} 1 &  \frac{A_{1j}}{A_{11}} \\ \\ 
0  & A_{ij} -\frac{A_{i1} A_{1j}}{A_{11}} \end{array}\right), 
$$
where $i \geq 2$ is the row index, $j\geq 2$ the column index
and $\delta_{ij}$ the Kronecker delta. 
It is immediate that ${\rm curl\,} M =0$, ${\rm div\,} P=0$, and, restoring
$\vep$ again, we have that $A^\vep \stackrel{H}{\longrightarrow} A$
if and only if $M^\vep \rightharpoonup M$, 
$P^\vep \rightharpoonup P$ weakly in $L^2 (\Om)^{n\times n}$.

%%%%%%%%%%%%%%%%%%%%%%%%%%%%%%%%%%%%%%%%%%%%%%%%%%%%%

\end{document}